\documentclass[12pt,twoside]{amsart}
\usepackage[all]{xy}
\usepackage{a4,amsmath,amssymb,amsfonts,amscd,mathrsfs}
\newtheorem{theorem}[subsection]{Theorem}
\newtheorem{lemma}[subsection]{Lemma}
\newtheorem{corollary}[subsection]{Corollary}
\newtheorem{proposition}[subsection]{Proposition}

\newcommand{\Db}{\mathrm{D^b}}
\newcommand{\DD}{\mathrm{D}}
\newcommand{\Ext}{{\rm Ext}}
\newcommand{\FM}{{\rm FM}}
\newcommand{\Hom}{{\rm Hom}}
\newcommand{\Jac}{{\rm Jac}}
\newcommand{\ol}[1]{\overline{#1}}
\newcommand{\Pic}{{\rm Pic}}
\newcommand{\Quot}{{\rm Quot}}
\newcommand{\SU}{{\rm SU}}
\newcommand{\U}{{\rm U}}
\newcommand{\ext}{{\rm ext}}
\newcommand{\im}{{\rm im}}
\newcommand{\pr}{{\rm pr}}
\newcommand{\rk}{{\rm rk}}
\newcommand{\Ecal}{{\mathcal E}}
\newcommand{\Hcal}{{\mathcal H}}
\newcommand{\Jcal}{{\mathcal J}}
\newcommand{\Ocal}{{\mathcal O}}
\newcommand{\Pcal}{{\mathcal P}}
\newcommand{\pdop}{{\mathbb P}}
\newcommand{\zdop}{{\mathbb Z}}
\newcommand{\dual}{^\lor}
\newcommand{\inv}{^{-1}}
\newcommand{\rarpa}[1]{\stackrel{#1}{\longrightarrow}}

\author{Georg Hein}
\begin{document}
\address{FB Mathematik, University Duisburg-Essen,
45117 Essen, Germany}
\email{georg.hein@uni-due.de}
\title{Minimal bundles and fine moduli spaces}
\date{March 13, 2009}
\maketitle
\begin{abstract}
We study sheaves $E$ on a smooth projective curve $X$ which are minimal
with respect to the property that $h^0(E \otimes L) >0$ for all line
bundles $L$ of degree zero. We show that these sheaves define ample
divisors $D(E)$ on the Picard torus $\Pic^0(X)$ (see Theorem
\ref{min-bun}). Next we classify all minimal sheaves of rank one (see
Theorem \ref{r1dg=min}) and two (see Theorem \ref{r2d-min-deg}).
As an application we show (see Proposition \ref{u=quot}) that the moduli
space parameterizing rank two bundles of odd degree can be obtained as a
Quot scheme.
\end{abstract}
\section{Introduction}
In this article we resume the study of minimal
bundles from \cite{Hei}. We call a sheaf $E$ on a
smooth projective curve $X$ minimal when $h^0(E \otimes L) >0$ for all
line bundles $L$ of degree zero, but $H^0(E' \otimes L)=0$ for all
proper subsheaves $E' \subsetneq E$ and a general line bundle $L$ of
degree zero. For more equivalent characterizations of minimality see
the definition after Theorem \ref{prop-Pm}.

Using the Fourier-Mukai transform and the Poincar\'e line bundle $\Pcal$
on $X \times \Pic^0(X)$ we find that for a minimal sheaf $E$ on $X$
the direct image sheaf $\pr_{2*}( \Pcal \otimes \pr_1^*E)$ is a line
bundle $\Ocal_\Pic(-D(E))$ where $D(E)$ is an effective ample divisor on
the Picard torus $\Pic^0(X)$. This is shown in Theorem \ref{min-bun}
which is basically contained in \cite[Theorem 4.6]{Hei}. However, the
proof given here (see Lemmas \ref{DE=ample} and \ref{zprime=0}) is
shorter than the old proof in \cite{Hei}.
In Proposition \ref{lotta-mins} we show that any sheaf $E$ on $X$ such
that the direct image sheaf $\pr_{2*} (\Pcal \otimes \pr_1^*E)$ is not
zero contains a minimal bundle. In Theorems \ref{r1dg=min} and
\ref{r2d-min-deg} we classify all minimal bundles of rank one and two.
As an unexpected application we obtain that the fine moduli schemes
$\Pic^d(X)$, $\U_X(2,2k+1)$ and $\SU_X(2,2k+1)$ are special cases of
Grothendieck's Quot schemes developed in \cite{Gro}. See Propositions 
\ref{pic=quot}, \ref{u=quot}, and \ref{su=quot} for the details.

The classical construction of geometric invariant theory (\cite{GIT},
see also \cite{New}) uses the Quot schemes as a starting point and
investigates a group action on this scheme. 
In the above cases we see that taking {\em the right Quot scheme} all
stable objects appear as quotients, all
quotients are stable, and for two quotients $[P \rarpa{\pi_1} E_1]$
and $[P \rarpa{\pi_2} E_2]$ we have an equivalence $E_1 \cong E_2 \iff
\ker(\pi_1)=\ker(\pi_2)$.

It may be well known that fine moduli spaces of vector bundles on
algebraic varieties which are projective appear as Quot schemes. For
lack of reference we include this result in Proposition \ref{fine=quot}.
In view of this result, the Propositions \ref{u=quot},
and \ref{su=quot} show that the Poincar\'e bundle on $X \times
\Pic^0(X)$ and the Quot scheme construction are enough to construct
$\U_X(2,2k+1)$ and $\SU_X(2,2k+1)$.

Since we exploit the properties of the Fourier-Mukai transform of the
Fourier-Mukai pair of abelian varieties $(\Jac(X),\Pic^0(X))$, Raynaud's
vector bundles $P_m$ (see Theorem \ref{prop-Pm}) play a crucial part.
These bundles are a special case of GV-sheaves defined by Pareschi and
Popa in \cite{PP}.

\subsection*{Acknowledgment}
This work has been supported by the SFB/TR 45
``Periods, moduli spaces and arithmetic of algebraic varieties''.

\section{Preliminaries}
\subsection*{Notation}
We fix an algebraically closed field $k$ of arbitrary characteristic.
In this article the word {\em scheme} means {\em scheme of finite type
over} $k$.
Let $X$ be a smooth projective curve of genus $g$ over $k$.
We fix a point $P_0 \in X(k)$.
We consider the following morphisms
\[\xymatrix{
& \Jac(X) \times \Pic^{0}(X) \ar[r]^-{p} \ar[d]_-{q} & \Pic^{0}(X)\\
X \ar@{^(->}[r]^-{\iota} & \Jac(X)\\
}\]
The Poincar\'e bundle on $\Jac(X) \times \Pic^0(X)$ is denoted by
$\Pcal$. We normalize it by demanding that
$\Pcal| _{ \{P_0 \} \times \Pic^0(X)} \cong \Ocal_{\Pic^0(X)}$. Mukai
showed in \cite{Muk} that the Fourier-Mukai transform 
\[ \FM_\Pcal : \Db(\Jac(X)) \to \Db(\Pic^0(X)) \quad e \mapsto
p_*(\Pcal\otimes q^*e) \]
is an equivalence of triangulated categories. Furthermore, he 
showed that the inverse $\FM\inv_\Pcal : \Db(\Pic^0(X)) \to \Db(\Jac(X))$
is given by $f \mapsto [-1]^*q_*(\Pcal \otimes p^*f)[g]$ where $[-1]$
denotes the involution on $\Jac(X)$ which assigns each element its
inverse. For a coherent sheaf $E$ on $X$ we may regard $\iota_*E$ as a
complex and obtain a complex $\FM_\Pcal(\iota_*E)$ in $\Db(\Pic^0(X))$.
For dimension reasons all the cohomology sheaves
$\FM^i(E):=\Hcal^i(\FM_\Pcal(\iota_*E))$ vanish for
$i \not \in \{ 0, 1 \}$.

\subsection*{The theta divisor $\Theta \subset \Pic^0(X)$}
We repeat some basic facts about the theta divisor on the Picard torus
(cf.~the first chapter of the book \cite{ACGH}).
We start with fixing a theta characteristic
$\vartheta$, that is a line bundle such that $\vartheta^{\otimes 2} \cong
\omega_X$.
This allows a set-theoretical definition of the $\Theta$-divisor
\[ \Theta = \{ [L] \in \Pic^0(X) \,|\, h^1(L \otimes \vartheta) >0 \} .\]
For a line bundle $M$ of degree $g-1$ we can likewise define
the theta divisor $\Theta_M$ by
$\Theta_M = \{ [L] \in \Pic^0(X) \,|\, h^1(L \otimes M) >0\}$.
$\Theta_M$ is the image of the summation morphism
from the $(g-1)$-fold symmetric product
$X^{(g-1)}$ of $X$ to $\Pic^0(X)$ given by
$(P_1,P_2,\dots,P_{g-1}) \mapsto \Ocal_X(P_1+P_2+\dots +P_{g-1}) \otimes
M\inv$.
$\Theta$ defines a principal polarization on $\Pic^0(X)$, that is
$\Theta$ is ample, and $\int_{\Pic^0(X)} \Theta^g = g!$.
The divisors $\Theta$ and $\Theta_M$ are not rationally equivalent
unless $M \cong \vartheta$. However, $\Theta_M$ is the translate of
$\Theta$ by $[\vartheta \otimes M\inv] \in \Pic^0(X)$.
Such divisors are called geometrically equivalent, we write $\Theta \sim
\Theta_M$.

\subsection*{Raynaud's bundles $P_m$}
For positive integers $m$ Raynaud considered in \cite{Ray} the vector
bundles $P_m = \iota^*[-1]^*R^gq_*(\Pcal \otimes p^*\Ocal_\Pic(-m\Theta))$.
We repeat some properties of $P_m$ next.
\begin{theorem}\label{prop-Pm}
Let $m$ be a positive integer, and $E$ a sheaf of rank $r$ on $X$.\\
\begin{tabular}{lp{13cm}}
(i) & $P_m$ is a vector bundle of rank $m^g$ and degree $gm^{g-1}$.\\
(ii) & There exists a unique (up to scalar multiplication) surjection
$P_{m+1} \to P_m$.\\
(iii) & We have equalities\\
&$\begin{array}{rcl}
\Hom_X(P_m,E) & =  & \Hom_{\Pic(X)}(\Ocal_\Pic(-m\Theta),\FM^0(E))\\
&=&\Hom_{\Db(\Pic(X))}(\Ocal_\Pic(-m\Theta),\FM_\Pcal(\iota_* E)).
\end{array}$\\
(iv) & The following equivalences hold\\
& $\begin{array}{rcl}
\FM^0(E) \ne 0
& \iff  & \Hom_{\Pic(X)}(\Ocal_\Pic(-m\Theta),\FM^0(E))
 \ne 0 \mbox{ for } m \gg 0\\
& \iff  & \Hom_{\Pic(X)}(\Ocal_\Pic(-m\Theta),\FM^0(E))
 \ne 0 \mbox{ for some } m > rg\\
& \iff & h^0(E \otimes L) >0 \mbox{ for all line bundles } L
\mbox{ of degree zero}\\
\end{array}$
\end{tabular}
\end{theorem}
\begin{proof}
(i) is 3.1 in Raynaud's article \cite{Ray}, (ii) is 2.4 in \cite{Hei},
(iii) is 2.5 in \cite{Hei}, and (iv) is the content of 2.5 and 3.7 in
\cite{Hei}. Partially, this is contained in Corollary 7.3 in \cite{PP}.
\end{proof}
\subsection*{Minimal bundles} We call a vector bundle $E$ on $X$
minimal, when $\FM^0(E) \ne 0$, and for all proper subsheaves $E'
\subsetneq E$ we have $\FM^0(E')=0$.
In view of Theorem \ref{prop-Pm} minimality can be expressed by the
following equivalent conditions\\
\begin{tabular}{rp{11.5cm}}
$E$ is minimal  $\iff$ & $\FM^0(E) \ne 0$, and $\FM^0(E')=0$ for all
proper subsheaves $E' \subset E$.\\
$\iff$ & $\Hom(P_m,E) \ne 0$, and $\Hom(P_m,E') =0$ for all proper
subsheaves $E' \subset E$, and $m \gg 0$.\\
$\iff$ & $\Hom(P_m,E) \ne 0$, and $\Hom(P_m,E') =0$ for all proper
subsheaves $E' \subset E$, and $m > rg$.\\
$\iff$ & $\Hom(P_m,E) \ne 0$ for some $m>rg$ and every $\pi \in
\Hom(P_m,E)$ is either zero or surjective.\\
$\iff$&$h^0(E \otimes L) \geq 1 $, and $h^0(E' \otimes L) = 0$ for all
proper subsheaves $E' \subset E$, and a general line bundle $L$ of degree
zero.\\
\end{tabular}
Note that the definition of {\em minimal bundle} differs slightly
from the definition given in \cite{Hei}. 

\section{$\FM^0(E)$ of minimal bundles}
Taking an elementary transformation
$0 \to E' \to E \to k(P) \to  0$,
we obtain an exact sequence
\[ 0 \to \FM^0(E') \to  \FM^0(E ) \to  L_P \]
where $L_P = \FM^0(k(P))$ is the numerical trivial line bundle
on $\Pic^0(X)$ parameterized by the point $P$.
The minimality of $E$ implies that $\FM^0(E')=0$, thus
\[ \FM^0(E) \subset L_P \mbox{ for all points } P \in X(k) .\]
In particular we have $\FM^0(E) \subset L_{P_0}=\Ocal_{\Pic^0(X)}$.
We fix an embedding $\FM^0(E) \to \Ocal_\Pic$, and
we write $\FM^0(E) = \Jcal_Z$ for a closed subscheme $Z \subset
\Pic^0(X)$.
We decompose $Z = D \cup Z'$, where $\Ocal_{D}$ is the quotient of
$\Ocal_Z$ modulo the maximal subsheaf of codimension two. This way we
can write $\Jcal_Z=\Ocal_\Pic(-D) \otimes \Jcal_{Z'}$ with $Z'$ of
codimension at least two.
\begin{lemma} \cite[Lemma 4.10]{Hei}\label{DE=ample}
The divisor $D$ is ample.
\end{lemma}
\begin{proof}
We consider the map of abelian varieties
$\Pic^0(X) \rarpa{\phi(D)} \Jac(X)$ given by
$\phi(D): P \mapsto T_P^* \Ocal(D) \otimes \Ocal(-D)$.
Its kernel $K(D)$ is the biggest subscheme where $D$ is trivial (see
\cite[Corollary 5 on page 131]{Mum}.)
If the kernel is discrete, then $D$ is ample, and we are done.
Let $K_0(D)$ be the connected component of zero in $K(D)$.
The injection $K_0(D) \to \Pic^0(X)$ defines a surjection 
$\Jac(X) \rarpa{\psi(D)} \Pic^0(K_0(D))$ of the dual abelian varieties.
If $L_P$ is the line bundle parameterized by the point $P \in X(k)$, then
we have $H^0(L_P \otimes \Ocal_\Pic(D)) = \Hom(\Ocal_\Pic(-D),L_P) =
\Hom(\Jcal_Z,L_P)$ because the codimension of $Z'$ is at least two.
As we remarked above there exists an embedding $\Jcal_Z = \FM^0(E)
\subset L_P$ for all points $P \in X(k)$. The global section of $H^0(L_P
\otimes \Ocal_\Pic(D))$ shows that $L_P \otimes \Ocal_\Pic(D)$ is
trivial on all fibers of $\phi(D)$. So it is trivial on $K_0(D)$.
We conclude that the kernel of $\psi(D)$ contains all points $P \in X(k)
\subset \Jac(X)(k)$. Since there exists no proper abelian subvariety in
$\Jac(X)$ which contains $X$, we conclude that the kernel of $\psi(D)$
is $\Jac(X)$. Thus, $K_0(D)$ is of dimension zero which implies the
ampleness of $D$.
\end{proof}

\begin{lemma}\cite[Lemma 4.11]{Hei}\label{zprime=0}
$Z'= \emptyset$ which implies $\FM^0(E)= \Ocal(-D)$.
\end{lemma}
\begin{proof}
We start with the short exact sequence $0 \to \Jcal_Z \to
\Ocal_\Pic(-D) \to \Ocal_{Z'} (-D) \to 0$.
Since $D$ is ample we have $R^iq_*(\Pcal \otimes p^* \Ocal_\Pic(-D))=0$
for all $i \ne g$. Since the dimension of $Z'$ is at most $g-2$ we have
$R^iq_*(\Pcal \otimes p^* \Ocal_{Z'}(-D))=0$ for $i \geq g-1$.
It follows that $\FM_\Pcal\inv(\Jcal_Z)$ has cohomology only in
nonpositive degrees, and the zero cohomology
$\Hcal^0(\FM_\Pcal\inv(\Jcal_Z)) =
\Hcal^0(\FM_\Pcal\inv(\Ocal_\Pic(-D)))$.
Since $\Jcal_Z=\FM^0(E)$ we have a canonical morphism
$\Jcal_Z \to \FM_\Pcal(E)$. We obtain a morphism
$\rho: \FM_\Pcal\inv \Jcal_Z \to \iota_*E$. The only non vanishing term in the
Eilenberg-Moore spectral sequence (see \cite[page 263]{GM})
\[\Ext^q(\Hcal^{q}(\FM_\Pcal\inv(\Jcal_Z)),\iota_*E)) \implies
\Hom(\FM_\Pcal\inv(\Jcal_Z),\iota_*E) = \Hom(\Jcal_Z,\FM_\Pcal(E))\]
is $\Hom(\Hcal^0(\FM_\Pcal\inv \Jcal_Z),\iota_*E)$. So we conclude that
$\rho$ factors through $\FM_\Pcal\inv(\Ocal_\Pic(-D))$. So $\Jcal_Z \to
\FM^0(E)$ factors through $\Ocal_\Pic(-D)$ which gives $\Jcal_Z =
\Ocal_\Pic(-D)$.
\end{proof}

\begin{theorem}\label{min-bun}
Let $E$ be a minimal bundle on the smooth projective curve $X$.\\
\begin{tabular}{lp{14cm}}
(i) & $\FM^0(E) = \Ocal_\Pic(-D)$ where $D = D(E)$ is an ample divisor.\\
(ii) & $\FM^1(E)$ has projective dimension at most two.\\
(iii) & $\FM\inv_\Pcal(\Ocal_\Pic(-D(E)))$ is a vector bundle of rank
$h^0(\Ocal_\Pic(D(E)))$, and the
natural map $\iota^* \FM\inv_\Pcal(\Ocal_\Pic(-D(E))) \to E$ is
surjective.\\
(iv) & $E$ is simple.\\
\end{tabular}
\end{theorem}
\begin{proof}
(i) follows from Lemma \ref{DE=ample} and \ref{zprime=0}.\\
In order to see (ii) consider a surjection
$M^{\oplus N} \rarpa \pi E$ with $M$
a line bundle of negative degree. By base change we obtain that
$\FM^1(M^{\oplus N} )$ and $\FM^1(\ker \pi)$ are vector bundles, and
$\FM^0(M^{\oplus N} ) = 0$. Hence, we obtain a resolution of $\FM^1(E)$
by vector bundles
\[0 \to \FM^0(E) \to \FM^1(\ker \pi) \to \FM^1(M^{\oplus N} ) \to
\FM^1(E) \to 0\]
(iii) By Serre duality
$\FM\inv_\Pcal(\Ocal_\Pic(-D(E))) = [-1]^*R^gq_*(\Pcal \otimes
p^*\Ocal_\Pic(-D(E)))$ is a vector bundles of rank
$h^0(\Ocal_\Pic(D(E)))$. Since $\Hom(P_m,E)= \Hom(\Ocal_\Pic(-m\Theta),
\FM^0(E))$, every morphism $P_m \to E$ factors through the image of
$\iota^* \FM\inv_\Pcal(\Ocal_\Pic(-D(E))) \to E$.\\
(iv) From the surjection $\iota^* \FM\inv_\Pcal(\Ocal_\Pic(-D(E))) \to
E$ we deduce that $\Hom(E,E) \subset \Hom(\iota^*
\FM\inv_\Pcal(\Ocal_\Pic(-D(E))),E)$. Since $\FM_\Pcal$ is an
equivalence of categories the last vector space is isomorphic to 
$\Hom(\Ocal_\Pic(-D(E)), \FM_\Pcal(\iota_*E)) = \Hom(\Ocal_\Pic(-D(E)),
\FM^0(E))$ which is one dimensional by definition of $D(E)$.
\end{proof}

\section{Examples of minimal bundles}
\begin{proposition}\label{lotta-mins}
Any coherent sheaf $E$ on $X$ with $\FM^0(E) \ne 0$ contains a minimal
sheaf.
\end{proposition}
\begin{proof}
Set $m=\rk(E) g +1$.
The non trivial homomorphisms $\phi: P_m \to E$ form a bounded family.
So there exists a $\phi_{\min}$ such that the image $E'=
\im(\phi_{\min})$ has minimal Hilbert polynomial. By construction any
morphism $P_m \to E'$ is either zero or surjective. In consequence, $E'$
is minimal.
\end{proof}

\begin{remark}
So any sheaf $E$ of positive Euler characteristic $\chi(E)$ contains a
minimal sheaf. This minimal sheaf is by no means uniquely defined.
As an example take a line bundle $M$ of degree $2g$. This line bundle
contains all line bundles $L$ of degree $g$ which are all minimal as we
show next.
\end{remark}

\begin{theorem}\label{r1dg=min}
A line bundle $L$ is minimal if and only if its degree $\deg(L)=g$.
For a line bundle $L$ of degree $g$ the divisor $D=D(L)$ is
geometrically equivalent to the $\Theta$-divisor on $\Pic^0(X)$.
\end{theorem}
\begin{proof}
If $\chi(L) \leq 0$, then a general line bundle $L'$ of degree $\deg(L)$
has no sections. Thus, we have $h^0( L \otimes (L' \otimes L\inv))=0$
which implies $\FM^0(L)=0$. If $\chi(L) >1$, then there are proper
subsheaves $L' \subset L$ of degree $g$ for which $\FM^0(L') \ne 0$.
For $\deg(L)=g$, we have seen that $\FM^0(L) \ne 0$, and all proper
subsheaves $L' \subset L$ satisfy $\FM^0(L')=0$. So the minimality is
shown.\\
Let now be $L$ a line bundle of degree $g$. From $0 \to L(-P_0) \to L
\to k(P_0) \to 0$ we deduce that 
\[ (\FM^0(L(-P_0)) =0) \to (\FM^0(L) = \Ocal_\Pic(-D(L))) \to \Ocal_\Pic
\to \FM^1(L(-P_0)) \,. \]
$\FM^1(L(-P_0))$ is a sheaf of projective dimension one concentrated on
the $\Theta$ divisor $\Theta_{L(-P_0)}$ which is geometrically
equivalent to $\Theta$. The linear system $|\Theta_{L(-P_0)}|$ consists
only of the divisor $\Theta_{L(-P_0)}$ which is a translate of the
irreducible $\Theta$ divisor. Since $\Theta_{L(-P_0)}$ is
irreducible and $D(L)$ is ample, we have $D(L)=\Theta_{L(-P_0)}$.
\end{proof}

\begin{theorem}\label{r2d-min-deg}
A rank two vector bundle $E$ is minimal if and only if it is stable and
of degree $2g-1$. For a minimal rank two bundle $E$ the divisor $D(E)$
is geometrically equivalent to $2\Theta$.
\end{theorem}
The statement of the proposition is a consequence of the three successive
Lemmas \ref{r2d=min}--\ref{r2d-deg}.
\begin{lemma}\label{r2d=min}
Any stable rank two vector bundle bundle $E$ of degree $2g-1$ is
minimal.
\end{lemma}
\begin{proof}
Suppose $E$ is a stable vector bundle of rank two with $\deg(E)=2g-1$.
Let $\pi:E \to k(P)$ be a surjection and let $E'$ be its kernel.
The stability of $E$ implies the semistability of $E'$. Therefore (see
\cite[Proposition 1.6.2]{Ray}) the sheaf $\FM^0(E')$ is zero.
Since $\chi(E) =1$ we have $\FM^0(E) \ne 0$. However any proper subsheaf
$E'' \subset E$ is contained in a sheaf $E'$ as above which implies
$\FM(E'')=0$. This gives the minimality of $E$.
\end{proof}

\begin{lemma}\label{r2min=stab}
Let $E$ be a minimal vector bundle of rank two.
Then $E$ is stable and of degree $\deg(E)=2g-1$.
\end{lemma}
\begin{proof}
First let us assume that $E$ were not stable. Then we would have a short
exact sequence $0 \to L_1 \to E \to L_2 \to 0$ with $\deg(L_1) \geq
\deg(L_2)$. Thus, the vanishing of $\FM^0(L_1)$ implies that
$\FM^0(L_2)=0$. Since $\FM^0(E) \ne 0$, we conclude that
$\FM^0(L_1) \ne 0$
which contradicts the minimality of $E$. We deduce that $E$ is stable.\\
Raynaud showed in Proposition 1.6.2 in \cite{Ray} that for a semistable
$E$ we have an equivalence $\FM^0(E) \ne 0 \iff \chi(E) > 0$. So we
deduce by minimality that $\chi(E)=1$. Riemann-Roch gives the stated
degree.
\end{proof}

\begin{lemma}\label{r2d-deg}
Let $E$ be a minimal bundle of rank two.
The divisor $D=D(E)$ is geometrically equivalent to the the
double of the $\Theta$-divisor on $\Pic^0(X)$.
\end{lemma}
\begin{proof}
By Lemma \ref{r2min=stab} the vector bundle $E$ is stable of degree
$\deg(E)=2g-1$. We consider a surjection $\pi: E \to k(P_0)$.
The bundle $E'=\ker(\pi)$ is semistable of Euler characteristic
$\chi(E')=0$. Thus, we have that $\FM^0(E')=0$ and
$c_1(\FM^1(E')) \sim 2 \Theta$.
Once we show that the support of $\FM^1(E)$ is at most $g-2$
dimensional we conclude, as in the proof of Theorem \ref{r1dg=min}
that $D(E) \sim 2\Theta$.\\
Case 1: $E$ contains a line bundle $L_1$ of degree $g-1$.\\
We choose the surjection $\pi$ in such a way that $\pi|_{L_1} =0$. This
way $L_1 \subset E'$. Therefore $E'$ is properly semistable, and the
$\Theta$-divisor of $E'$ is the union of the irreducible divisors
$\Theta_{L_1}$ and $\Theta_{E/L_1}$. From the short exact sequence
$0 \to E' \to E \to k(P_0) \to 0$ we deduce the long exact sequence
\[ 0 \to \left( \FM^0(E) = \Ocal_\Pic(-D(E)) \right) \to \Ocal_\Pic \to
\FM^1(E') \to \FM^1(E) .\]
Since the support of $\FM^1(E')$ is the union of $\Theta_{L_1}$ and
$\Theta_{E/L_1}$ we conclude that $D(E)$ is one of the four divisors
in the set $\{ \emptyset, \Theta_{L_1}, \Theta_{E/L_1}, \Theta_{L_1}
+ \Theta_{E/L_1} \}$. Since $D(E)$ is ample (see Theorem
\ref{min-bun} (i)) we can exclude $\emptyset$. We can exclude
$\Theta_{L_1}$ because $h^0(\Ocal_\Pic(\Theta_{L_1}))=1$. So there can
not be a surjection from $\iota^*\FM\inv_\Pcal(\Ocal_\Pic(\Theta_{L_1}))$
to $E$ (see Theorem \ref{min-bun} (iii)). For the same reason we can
exclude $\Theta_{E/L_1}$. So we have $D(E)= \Theta_{L_1} +
\Theta_{E/L_1}$.

Case 2: $E$ contains no line bundle of degree $g-1$.\\
A line bundle $[L] \in \Pic^0(X)(k)$ is in the support of
$\FM^1(E) \iff H^1(E
\otimes L) \ne 0$. By Serre duality this is equivalent to
$\Hom(E,\omega_X \otimes L\inv) \ne 0$. So it is enough to show that the
line bundles $M$ of degree $2g-2$ which admit a morphism
$E \rarpa{ \ne 0} M$ form a
family of dimension at most $g-2$. Now any morphism $E \rarpa{\alpha} M$
factors through its image $E \to \im(\alpha) \to M$.
Let $\deg(\im(\alpha)) = 2g-2-d$. The semistability of $E$ implies 
$d \leq g-2$. Since for a given line bundle $\im(\alpha)$ the extensions 
$\im(\alpha) \to M$ are parameterized by effective divisors of degree
$d$ we have to show that the image of the Quot scheme
$\Quot^{1,2g-2-d}_E$ of line bundle quotients of $E$ of degree $2g-2-d$
is of dimension at most $g-2-d$.
We consider the following morphisms
\[ \xymatrix{ \Quot^{1,2g-2-d}_E \ar[rrr]_-{[E \to L] \mapsto L}^-\beta
&&& \Pic^{2g-2-d}(X)
\ar[rrr]_-{L \mapsto L^{\otimes 2} \otimes \det(E)\inv}^-\gamma
&&& \Pic^{2g-3-2d}(X) } \]
Since $\gamma$ is a finite morphism we have to show that the image of
$\gamma \circ \beta$ is of dimension at most $g-2-d$.

The Quot scheme parameterizing quotients $[E \rarpa{\pi} L]$ has
tangent space
$\Hom(\ker(\pi),L)=H^0(L\otimes \ker(\pi)\inv)$.
The line bundle
$L\otimes \ker(\pi)\inv
 \cong L^{\otimes 2} \otimes \det(E)\inv$ is of degree $2g-3-2d$.
From Clifford's theorem \cite[Theorem IV.5.5]{Har} we see that 
the $h^0(L^{\otimes 2} \otimes \det(E)\inv) \leq g-1-d$.
We assume now that there exists a irreducible component $Q \subset
\Quot^{1,2g-2-d}_E$ of dimension $g-1-d$ such that its image 
$\ol Q =\gamma \circ \beta(Q)$ is also of dimension $g-1-d$.
We consider the cartesian diagram
\[ \xymatrix{ Z \ar@{=}[r]
& \ol{Q} \times _\Pic X^{(2g-3-2d)} \ar[r] \ar[d]_\xi
& X^{(2g-3-2d)} \ar[d]\\
& \ol Q \ar[r] & \Pic^{2g-3-2d} }\]
Since for all line bundles $\ol L$ parameterized by $\ol Q$ we have
$h^0(L)= g-1-d$ we follow that $\xi$ is a $\pdop^{g-2-d}$ bundle. So $Z$
is of dimension $2g-3-d$. So we must have $Z= X^{(2g-3-2d)}$. However by
Poincar\'e's formula (cf I.\S 5 in \cite{ACGH})
the image of $X^{(2g-3-2d)}$ is of dimension $2g-3-2d$.
So we conclude that $g-2-d=0$ which implies that $d=g-2$. This means
that $E$ has a quotient line bundle $[E \rarpa \pi L]$ of degree
$\deg(L)=g$. The kernel of $\pi$ is a line bundle of degree $g-1$ which
we excluded  in case 2.
\end{proof}

\section{Quot Schemes as fine Moduli Spaces}
\subsection{The general case}
Let $(X,\Ocal_X(1))$ be a polarized projective variety, and $\chi:\zdop
\to \zdop$ a polynomial. We say that a scheme $M$ is up to first order
a fine moduli space
of simple $\chi$-bundles on $X$ when the following conditions are
satisfied:
\begin{tabular}{lp{14cm}}
(i) & There exists a vector bundle $\Ecal$ on $X \times M$.\\
(ii) & For any point $P \in M(k)$ the vector bundle $\Ecal_P:= \pr_{1*}(
\Ecal \otimes \pr_2^*k(P))$ has Hilbert polynomial $\chi$ with respect to
$\Ocal_X(1)$.\\
(iii) & For $P,Q \in M(k)$ we have $\hom(\Ecal_P, \Ecal_Q) =
\left\{ \begin{array}{ll} 1 & \mbox{when } P=Q\\ 0 & \mbox{otherwise.}
\end{array} \right.$\\
(iv)& For any $P \in M(k)$ the Kodaira-Spencer map $T_{M,P} \to
\Ext^1(\Ecal_P,\Ecal_P)$ is an isomorphism.
\end{tabular}

Note that any open subset of a fine moduli space is a up to first order
a fine moduli space.
The following notion is very convenient: For a point $P \in M(k)$ we
have the vector bundle $E = \Ecal_P$ parameterized by $P$ on $X$.
Therefore we denote the point $P$ by $[E]$.  

\begin{proposition}\label{fine=quot}
Let $(X,\Ocal_X(1))$ be a polarized projective Gorenstein variety, and $\chi:\zdop
\to \zdop$ a polynomial.
Assume there exists a up to first order fine moduli space  $M=M_X^\chi$ of simple
$\chi$-bundles on $X$. Furthermore, we assume that
$M$ is smooth, connected, and projective over $k$.
Under these assumptions there exists a vector bundle $P$ on $X$ such that
a connected component $Q_1$ of the Quot scheme $\Quot_{P/X}^\chi$ of
quotients of $P$ with Hilbert polynomial $\chi$ gives an isomorphism
\[ \xymatrix{ Q_1 \ar[rrr]^-\sim_-{ [Q \to E] \mapsto [E]} &&& M }\]
which identifies the universal quotient $\Ecal$ of $\pr_1^*P$ on $X
\times Q_1$ with the universal bundle on $X \times M$.
\end{proposition}
\begin{proof}
First, we fix an ample polarization $\Ocal_M(1)$ on $M$.
We consider the morphisms
$\xymatrix{X & X \times M \ar[r]^-q \ar[l]_-p & M}$. Any twist
$\Ecal(N) := \Ecal \otimes q^*\Ocal_M(N)$ is also a universal bundle on
$X \times M$. For a closed point $[E] \in M(k)$ with ideal sheaf
$\Jcal_E \subset \Ocal_M$ we obtain a long exact sequence on $X$:
\[ \dots \to p_*\Ecal(N) \to p_*(\Ecal(N) \otimes q^*k([E])) \to 
R^1p_* (\Ecal(N) \otimes q^*\Jcal_E) \to \dots \]
We may choose $N \gg 0$ such that the higher direct image sheaves $R^ip_*
(\Ecal(N) \otimes q^*\Jcal_E)$ all vanish for all points $[E] \in M(k)$,
and all $i \geq 1$.
This forces $p_*\Ecal(N)$ to be a vector bundle. By definition the
vector bundle $E$ is isomorphic to $p_*(\Ecal(N) \otimes q^*k([E]))$.
Summing up, we obtain a surjection $p_*\Ecal(N) \twoheadrightarrow
E$ for all $[E] \in M(k)$.\\
Next we want to show that $\Hom(p_*\Ecal(N),E)$ is one dimensional
given that $N$ is sufficiently big.
By Serre duality $\Hom(p_*\Ecal(N),E)$ is dual to
$\Ext^{\dim X}(E \otimes \omega_X \dual, p_*\Ecal(N))$.
Since $p^*$ is the left adjoint of $p_*$ we find
$\Ext^{\dim X}(E \otimes \omega_X \dual, p_*\Ecal(N)) =
\Ext^{\dim X}(p^*(E \otimes \omega_X \dual),\Ecal(N))$.
This vector space is isomorphic to
$H^{\dim X}(p^*E \dual \otimes p^*\omega_X \otimes \Ecal(N))$.
To compute this cohomology group we use the Leray spectral sequence
$$ H^k(R^lq_*(p^*E \dual \otimes p^*\omega_X \otimes \Ecal(N))) \implies
H^{k+l}(p^*E \dual \otimes p^*\omega_X \otimes \Ecal(N)) .$$
Choosing $N$ sufficiently high, all the higher cohomology groups $H^k( \dots
) $ on the left vanish and we obtain
that $H^{\dim X}(p^*E \dual \otimes p^*\omega_X \otimes \Ecal(N)) =
H^0(R^{\dim X}q_*(p^*E \dual \otimes p^*\omega_X \otimes \Ecal(N)))$.
By Serre duality, base change, and our condition $(iii)$ we see that $R^{\dim
X}q_*(p^*E \dual \otimes p^*\omega_X \otimes \Ecal(N))$ is the
skyscraper sheaf $k([E])$. This proves that
$\hom(p_*(\Ecal(N)),E)=1$.\\
Since the vector bundles parameterized by $M$ form a bounded family, we
may choose $N$ such that $\hom(p_*(\Ecal(N)),E)=1$ for all $[E ] \in
M(k)$.\\
Next we construct a morphism from $M$ to the Quot scheme
$\Quot^\chi_{p_*\Ecal(N)/X}$ as follows:
On $X \times M \times M$ we consider the surjection
$\pr_{13}^*\Ecal(N) \to \pr_{13}^*\Ecal(N) \otimes \pr_{23}
\Ocal_{\Delta(M)}$. Applying $\pr_{12*}$ to this surjection we obtain a
surjection $\pr_{12*} \pr_{13}^*\Ecal(N) \twoheadrightarrow \pr_{12*}
(\pr_{13}^*\Ecal(N) \otimes \pr_{23} \Ocal_{\Delta(M)})$.
The direct image sheaf $\pr_{12*} \pr_{13}^*\Ecal(N)$ is isomorphic to
$p^* p_* \Ecal(N)$, and the quotient $\pr_{12*}
(\pr_{13}^*\Ecal(N) \otimes \pr_{23} \Ocal_{\Delta(M)})$ becomes
isomorphic to $E$ when specialized to $X \times \{[E]\}$.
This way we obtain a morphism from $M \to \Quot^\chi_{p_*\Ecal(N)/X}$
which sends $[E]$ to $[p_*\Ecal(N) \rarpa \pi E]$ where $\pi$ is the
unique (up to scalar) surjective morphism in $\Hom(p_*\Ecal(N),E)$.\\
To finish the proof we show that the constructed morphism $M \to
\Quot^\chi_{p_*\Ecal(N)/X}$ identifies $M$ with one connected component
of $\Quot^\chi_{p_*\Ecal(N)/X}$. By replacing the Quot scheme by the
connected component containing the image of $M$, we may assume
$\Quot^\chi_{p_*\Ecal(N)/X}$ is connected. Since all the sheaves
parameterized by $M$ are different the morphism $M \to
\Quot^\chi_{p_*\Ecal(N)/X}$ is injective on closed points. Thus, $\dim
\Quot^\chi_{p_*\Ecal(N)/X} \geq \dim M$. 
For any $[E] \in M(k)$ we consider the surjection $p_*\Ecal(N) \rarpa
\pi E$. This gives a long exact sequence
\[ 0 \to \Hom(E,E) \rarpa \alpha \Hom(p_*\Ecal(N),E) \rarpa \beta
\Hom(\ker \pi, E ) \rarpa {{\rm d}_E} \Ext^1(E,E) \to \dots  \]
The morphism $\alpha$ is a morphism between two vector spaces of
dimension one. We deduce that $\beta$ is zero, and ${\rm d}_E$ is
injective. ${\rm d}_E$ is the tangent map for the forgetful morphism
$\Quot^\chi_{p_*\Ecal(N)/X} \to M$ which assigns $[p_*\Ecal(N)
\twoheadrightarrow E] \mapsto [E]$ in the point $[p_*\Ecal(N)
\twoheadrightarrow E]$.
We obtain
\[ \dim \Quot^\chi_{p_*\Ecal(N)/X} \leq \hom(\ker \pi, E ) \leq
\ext^1(E,E) = \dim(M) . \]
Putting together this inequality with the one obtained before ($\dim
\Quot^\chi_{p_*\Ecal(N)/X} \geq \dim M$) we see that all $\leq $ become
equalities and ${\rm d}_E$ is an isomorphism. Since $M$ is projective
its image on $\Quot^\chi_{p_*\Ecal(N)/X}$ is projective.
We conclude that $M \to \Quot^\chi_{p_*\Ecal(N)/X}$ is an isomorphism.
\end{proof}

\begin{remark}
The assumptions of Proposition \ref{fine=quot} are fulfilled when $X$ is
a smooth projective curve and $M=\U_X(r,d)$ is the moduli space of rank
$r$ stable vector bundles of degree $d$ on $X$ provided that $r$ and $d$
are coprime integers. Therefore, we find (under the assumption
$(r,d)=1$) that all the moduli spaces
$\U_X(r,d)$ are quot schemes $\Quot^{r,d}_{P/X}$ for certain sheaves $P$
on $X$.
\end{remark}

\subsection{The case of algebraic curves}
From now on $X$ is a smooth projective curve over $k$ of genus $g$.
We show, how the minimality of certain bundles allows the construction
of the moduli spaces of these bundles. Whereas in Proposition
\ref{fine=quot} we can not say much about the vector bundle $P$, we see that
Raynaud's bundles $P_m$ can be taken in the situations of
\ref{pic=quot}--\ref{su=quot}.

\begin{proposition}\label{pic=quot}
The Quot scheme $\Quot= \Quot_{P_2/X}^{1,g}$ of quotients of the Raynaud
bundle $P_2$ of rank one and degree $d$ is isomorphic to $\Pic^g(X)$.
The universal quotient of $\pr_1^*P_2$ on $X \times \Quot$ is a
Poincar\'e bundle.
\end{proposition}
\begin{proof}
Let $\xymatrix{P_2 \ar@{->>}[r] &Q}$ be surjection to a sheaf $Q$ on $X$
with $\rk(Q)=1$, and $\deg(Q)=g$. If $Q$ has torsion subsheaf $T$
different from zero, then there exists a surjection
$\xymatrix{P_2 \ar@{->>}[r] & L=Q/T}$
to a line bundle of degree $\deg(L) \leq g-1$. However, $L$ is properly
contained in a line bundle $M$ of degree $g$ which contradicts
Theorem \ref{r1dg=min}. We conclude that $Q$ is torsion free, hence
a line bundle.\\
Let now $Q$ be a line bundle of degree $g$ on $X$. We have seen that
$\FM^0(Q)= \Ocal_\Pic(-D)$ for some divisor $D \sim \Theta$. We
conclude that $\Hom(\Ocal_\Pic(-2\Theta),\FM^0(E))$ is one dimensional.
This implies (see (iii) of Theorem \ref{prop-Pm}) that $\Hom(P_2,E)$ is
one dimensional. Thus, we have a nontrivial morphism $\pi:P_2 \to Q$.
The minimality of $Q$ (Theorem \ref{r1dg=min}) implies that the
image of $\pi$ cannot be contained in a proper subsheaf. We eventually
conclude that $\Hom(P_2,Q)$ is a one dimensional vector space, generated
by an element $\pi$ which is surjective.\\ Thus, every point
$[\xymatrix{P_2 \ar@{->>}[r]^\pi & Q}]$ of $\Quot(k)$ corresponds to
exactly one line bundle $[Q] \in \Pic^d(X)(k)$. Let $\Quot=\Quot_1 \cup
\dots \cup \Quot_M$ be the decomposition into irreducible components.
Since the forgetful morphism $[\xymatrix{P_2 \ar@{->>}[r]^\pi & Q}]
\mapsto [Q]$ is surjective, we may assume that $\dim \Quot_1 \geq g =
\dim \Pic^d(X)$. Take a closed point $[\xymatrix{P_2 \ar@{->>}[r]^\pi &
Q}] \in \Quot_1(k)$. We have a long exact sequence
\[ 0 \to \Hom(Q,Q) \rarpa \alpha \Hom(P_2,Q) \rarpa \beta
\Hom(\ker(\pi),Q) \rarpa \gamma \Ext^1(Q,Q) \, .\]
We have seen that
$\Hom(P_2,Q)$ is one dimensional, so it follows that
$\alpha$ is an isomorphism, $\beta=0$, and $\gamma$ is injective.
$\Hom(\ker(\pi),Q)$ is the tangent space of $\Quot_1$ at $[\xymatrix{P_2
\ar@{->>}[r]^\pi & Q}]$, so it is of dimension at least $g$.
$\Ext^1(Q,Q)$ is the tangent space of $\Pic^d(X)$ at the point $[Q]$, a
$g$-dimensional vector space. Thus, we see that $\hom(\ker \pi , Q)=g$
for all points in $\Quot_1$, and $\gamma$ is an isomorphism. This
implies that $\Quot_1$ is regular.
We see that $\Quot_1 \to \Pic^d(X)$ is injective on closed points and an
isomorphism on the tangent spaces. Since $\Quot_1$ is projective, we
deduce that $\Quot = \Quot_1 \rarpa{\sim} \Pic^d(X)$.
\end{proof}

\begin{corollary}
For any integer $d \in \zdop$ there exists a sheaf $P^d_2$ such that the
Quot scheme $\Quot_d= \Quot^{1,d}_{P^d_2/X}$ of rank one quotients of
$P^d_2$ is isomorphic to $\Pic^d(X)$ and the universal quotient of
$\pr^*_1 P^d_2$ on $X \times \Quot_d$ is a Poincar\'e bundle.
\end{corollary}
\begin{proof} Indeed let $L$ be any line bundle of degree one, set
$P^d_2= P_2 \otimes L^ {\otimes d-g}$,
then the result follows from Proposition \ref{pic=quot}.
\end{proof}

\begin{proposition}\label{u=quot}
The Quot scheme $\Quot= \Quot^{2,2g-1}_{P_3/X}$ of rank two quotients of
 degree $2g-1$ of the Raynaud bundle $P_3$ is isomorphic to the moduli
space $\U_X(2,2g-1)$ of stable rank two bundles of degree $2g-1$. The
universal quotient of $\pr_1^*P_3$ on $X \times \Quot$ is a Poincar\'e
bundle.
\end{proposition}
\begin{proof}
Let $\xymatrix{P_3 \ar@{->>}[r]^\pi & F}$ be a surjection to a sheaf of
rank two of degree $2g-1$. Let $T$ be the torsion part of $F$, and set
$F''=F/T$.\\
Step 1: $F''$ does not surject to a line bundle $M$
of degree smaller than $g$.\\
If there were such a surjection $F'' \to M$,
then we could embed $M$ as a proper subsheaf of a line bundle $M'$ of
degree $g$, and would obtain from the composition
$\xymatrix{P_3 \ar@{->>}[r] & F \ar@{->>}[r] & F'' \ar@{->>}[r] & M
\ar[r]^{\subsetneq} & M'}$
a nonzero morphism $P_3 \to M'$ which is not surjective,
which contradicts Theorem \ref{r1dg=min}.\\
Step 2: We take now an elementary transformation $F'' \to F'$ such that
$\deg(F') = 2g-1$.\\ 
Step 3: $F'$ is stable.\\
Indeed, suppose that $F' \to M$ would be a destabilizing quotient. Then
$\deg(M) = \mu(M) \leq g-1$. The image of $F''$ in $M$ were a line
bundle of degree at most $g-1$ which is impossible by step 1.\\
Step 4: $T = 0$\\
Indeed if $T \ne 0$, then $P_3 \to F'$ would be a nontrivial morphism
which is not surjective which contradicts the minimality of $F'$ by
Theorem \ref{r2d-min-deg}.\\
We eventually conclude that $F = F'$ is stable.
Thus, we have seen that any quotient $[P_3 \to F]$ yields a stable
vector bundle $F$. Stable vector bundles $F$ (with $\rk(F)=2$, and
$\deg(F)=2g-1$) are by Theorem \ref{r2d-min-deg} quotients of $P_3$.
The rest follows along the lines of the proof of
Proposition \ref{pic=quot}.
\end{proof}

\begin{corollary}
For any odd integer $d$ there exists a vector bundle $P_3^d$
on $X$ such that 
the Quot scheme $\Quot= \Quot^{2,d}_{P_3^d/X}$ of rank two quotients of
 degree $d$ of the vector bundle $P_3^d$ is isomorphic to the moduli
space $\U_X(2,d)$ of stable rank two bundles of degree $d$. The
universal quotient of $\pr_1^*P^d_3$ on $X \times \Quot$ is a Poincar\'e
bundle.
\qed
\end{corollary}

\begin{proposition}\label{su=quot}
The Quot scheme $\Quot= \Quot^{2,2g-1}_{P_2/X}$ of rank two quotients of
 degree $2g-1$ of the Raynaud bundle $P_2$ is isomorphic to the moduli
space $\SU_X(2,\omega_X(P_0))$ of stable rank two bundles of determinant
isomorphic to $\omega_X(P_0)$. As before, the
universal quotient of $\pr_1^*P_2$ on $X \times \Quot$ is a Poincar\'e
bundle.
\end{proposition}
\begin{proof}
The birational divisor class $D(E)$ does not vary when $\det(E)$ stays
constant. Thus, all vector bundles $E$ of determinant $\omega_X(P_0)$
have the same divisor $D(E)$.
\end{proof}

\end{document}